\documentclass{amsart}
\usepackage{amssymb,amsmath,amsthm,amscd,float, amsaddr}
\usepackage{mathtools}

\usepackage{booktabs, topcapt}
\usepackage{enumerate}

\usepackage{hyperref}
\hypersetup{
	colorlinks   = true, 
	urlcolor     = blue, 
	linkcolor    = blue, 
	citecolor   = red 
}

\usepackage[sorted, msc-links]{amsrefs}

\DeclareMathOperator\wh{\mathfrak{h}}

\usepackage{cleveref}

\newtheorem{lem}{Lemma}

\newtheorem{rem}{Remark}
\newtheorem{exa}{Example}

\theoremstyle{definition}

\theoremstyle{remark}

\def\Z{\mathbb{Z}}
\def\Q{\mathbb{Q}}
\def\P{\mathbb{P}}
\def\R{\mathbb{R}}
\newcommand{\Ok}{\mathcal{O}_K}

\def\l{\lambda}
\def\X{\mathcal X}
\def\M{{\mathcal M}}
\def\p{\mathfrak{p}}
\def\h{\mathfrak{h}}
\def\w{\mathbf{w}}

\def\WP{\mathbb WP_{(2, 4, 6, 10)}^3}

\newcommand\A{\mathbb A}

\def\<{\langle}
\def\>{\rangle}

\def\awh{\mathfrak h}

\def\a{\alpha}

\def\l{\lambda}
\def\w{\mathfrak w}

\DeclareMathOperator\wgcd{wgcd }

\begin{document}

\title{The weighted moduli spaces of sextics}
	
\author{Lubjana Beshaj}
\address{Department of Mathematical Sciences, \\
 United States Military Academy, \\
 West Point, NY, 10996}
\email{lubjana.beshaj@usma.edu}

\author{Scott Guest}
\address{Department of Mathematics and Statistics,\\
 Oakland University, \\
 Rochester, MI 48309}
\email{sguest@oakland.edu}
	
\begin{abstract}
We use the weighted moduli height as defined in \cite{sh-h} to investigate the distribution of fine moduli points in the moduli space of genus two curves.
 We show that for any genus two curve with equation $y^2=f(x)$,  its weighted moduli height $\awh   (\p) \leq 2^3 \sqrt{3 \cdot 5 \cdot 7} \, \cdot  H(f)$, where $H(f)$ is the minimal naive height of the curve as defined in \cite{height}.
Based on the weighted moduli height $\awh$ we create a database of genus two curves defined over $\Q$ with small   $\awh$ and show that for small such height  ($\awh < 5$) about 30\% of points are fine moduli points.
\end{abstract}
	
	
\maketitle


\section{Introduction}

In \cite{height} authors studied heights of algebraic curves and compared the naive height with the moduli height.  Their examples focused mostly on curves of genus two.
In  \cite{data}, the authors created a database of rational points in the moduli space of genus two curves $\M_2$  with moduli height $\mathfrak h \leq 20$ and considered the problem of what percentage of these points have corresponding genus 2 curves defined over $\Q$.    It is widely believed that for  almost all points $\p \in \M_2 (\Q)$, there is no genus two representative defined over $\Q$. In the language of arithmetic geometry this means that for almost all genus two curves the field of moduli is not a minimal field of definition.
However,  results in \cite{data} were somewhat surprising. For all $\mathfrak h \leq 20$ roughly 30\% of the points the field of moduli was a field of definition. There were a couple interpretations to this.

Firstly, the points in  $\M_2 (\Q)$ were represented using the absolute invariants $i_1, i_2, i_3$ of genus two curves.   These are rational functions in terms of the Igusa invariants $J_2$, $J_4$, $J_6$, and $J_{10}$, see \cite{Ig-60}, \cite{data}, \cite{m-sh-2}, but they are not defined when $J_2 = 0$.  To handle the case $J_2 =0$ another compactification of the moduli space was used.  When using two different sets of invariants however, there is no natural means to make statements about interesting questions regarding the moduli space of genus two curves.

Secondly, it is possible that the "majority" of moduli points for which the field of moduli is a field of definition are concentrated around the "center" of the moduli space.  In other words, as the moduli height $\mathfrak h$ goes to $\infty$, the percentage of curves for which the field of moduli is a field of definition tends to zero.  However, computations in \cite{data}  didn't quite show this trend.

In \cite{m-sh} a "universal" equation of genus two curves was discovered.  Here universal means that it works for any tuple of arithmetic invariants $(J_2, J_4, J_6, J_{10})$ including the case when $J_2=0$ or the case when the automorphism group has order $> 2$, which was treated separately in Mestre's original algorithm \cite{Me} and has not been implemented in most computational algebra packages. This raised the question on whether we could perform computations done in \cite{data}, but now using $[J_2, J_4, J_6, J_{10}]$ and the weighted projective space $ \mathbb{WP}_{(2,4,6,10)}^3$ instead of the coordinate $(i_1, i_2, i_3)$ used in \cite{data}?

The main question became whether there was any measure of height for the weighted projective spaces which could be used for $ \mathbb{WP}_{(2,4,6,10)}^3$ in order to perform a similar analysis as in \cite{data}. Luckily, \cite{sh-h} came out which does precisely what we needed, defines a height on weighted projective spaces. Using this height for $ \mathbb{WP}_{(2,4,6,10)}^3$ we were able to determine a unique tuple $(J_2, J_4, J_6, J_{10})$ for any point $\p \in  \mathbb{WP}_{(2,4,6,10)}^3$ and compile a database of genus two curves based on this weighted projective height.

In this database every point is represented by a unique minimal tuple of arithmetic invariants $(J_2, J_4, J_6, J_{10})$.  Such database is valid for curves over any field of characteristic different from two.  In the case of characteristic two we can do exactly the same thing by adding another invariant $J_8$ as explained in \cite{Ig-60}.
The disadvantage of using the weighted moduli height is that since powers of invariants are involved, one can not perform computations for very big heights. In this paper we only go for heights up to $\mathfrak h =4$.

This paper is organized as follows.  In \cref{sec:3}, we start with giving a brief definition of a weighted projective space and then  define the height  function on this space, following closely the approach in  \cite{sh-h}.
	
In	\cref{sec:5}  we explain how we build the database. We start with creating a database of normalized tuples $\p = [J_2: J_4: J_6: J_{10}] \in \mathbb{WP}_{(2,4,6,10)}^3(\Q)$ with weighted height $\h \leq 5$ and $J_{10} \neq 0$.  Then,  in subsection	\ref{sec:4} we explain how  for each point $\p \in \M_2 (\Q) $ we can compute the equation of the curve over the rationals, when such exists. In \cite{m-sh} the authors do exactly this and construct for every point $\p \in \M_2$ the equation of the corresponding genus 2 curve. While we know that each $\p$ is defined over $\Q$ if it corresponds to a curve defined over $\Q$, the converse in general is not true. It is then necessary to determine if a given moduli point actually corresponds to a curve over $\Q$.    We say a moduli point $\p\in \M_2 (\Q)$ is a fine  point if and only we can find a representative curve $\mathcal \X_{\p}$ defined over $\Q$, otherwise we call it a coarse point or an obstruction point.
	
For each such tuple $\p = [J_2: J_4: J_6: J_{10}] \in \mathbb{WP}_{(2,4,6,10)}^3(\Q) \setminus  \{ J_{10} = 0\}$ , we also compute the automorphism group of the corresponding curve, using  algorithms presented in \cite{data}. We then present the  data on the distribution of fine points in the moduli space of genus two curves. Namely, we compare the number of moduli points below a given weighted moduli height to the number of fine points below that height for all weighted moduli heights $\h < 5$.  Finally, we give a table of specific results for all moduli points with $\h < 1$.

\bigskip

\noindent \textbf{Notation:}  Throughout this paper by a curve we mean a smooth, irreducible algebraic curve. Unless otherwise noted a curve $\X$ means the isomorphism class of $\X$ over some field $k$. With $ \mathbb{WP}_{(2, 4, 6, 10)}^3 (\Q)$  we denote the weighted projective space with weights $(2, 4, 6, 10)$ and  $\wh (\p)$ we denote  the weighted moduli height of a point $\p \in \  \mathbb{WP}_{(2, 4, 6, 10)}^3 (\Q)$ which we define in the preliminaries.

\section{Preliminaries}\label{sec:2}

Let $\X$ be a genus 2 curve defined over a field $k$ with characteristic 0 or a prime $p \neq 2$. It is a well known fact that the  equation of $\X$ can  be given by
\begin{equation}\label{eq-super}
\begin{split}
\X:  y^2 &=  f(x) = a_6 \sum_{i=1}^6 (x - \a_i) =   a_6 x^6 +  a_5 x^5 z+ \cdots + a_0 z^6 \\
\end{split}\end{equation}
where $\{\a_i\}_{i=1}^6$ are the ramification points of the map $\X \to \P^1$.   The isomorphism classes of genus 2 curves are on one to one correspondence with the orbits of the $GL_2 (k)$-action on the space of binary sextics.  Igusa proved that the invariant ring $\R_6$ is generated by the Igusa arithmetic invariants $J_2, J_4, J_6, J_{10}$; see   \cite{Ig-60, data} for details.  Note that Igusa $J$-invariants $\{J_{2i}\}$, for  $i=1, 2, 3, 5$, are homogenous polynomials of degree $2i$ in   $k[a_0, \dots, a_6]$ and $J_{10}$ is the discriminant of the curve.  Since in this paper we are working with smooth, irreducible curves than we take $J_{10} \neq 0$.  The definitions of these invariants can be found in several papers such as  \cite{Ig-60, data, deg2}.  %
%
%
\subsection{The weighted projective space}\label{sec:3}

In this subsection we will follow closely  the notation in \cite{sh-h}.     Let $k$ be a field of characteristic zero and  $(q_0, \dots , q_n)$ a fixed tuple of positive integers called \textbf{weights}.  Consider the action of $k^\star$ on $\A^{n+1}$ as follows
\[ \lambda * (x_0, \dots , x_n) = \left( \l^{q_0} x_0, \dots , \l^{q_n} x_n   \right) \]
for $\l\in k^\ast$.
The quotient of this action is called a \textbf{weighted projective space} and denoted by   $\mathbb{WP}^n_{(q_0, \dots , q_n)}$.
It is the projective variety $Proj \left( k [x_0,...,x_n] \right)$ associated to the graded ring $k [x_0, \dots ,x_n]$ where the variable $x_i$ has degree $q_i$ for $i=0, \dots , n$.

In \cite{sh-h} they define the concept of height on the weighted projective spaces. In this subsection we will mostly focus on the space $\mathbb{WP}_{(2,4,6,10)}^3(\Ok)\setminus \{J_{10} = 0 \}$.  Let  $\p = [J_2, J_4, J_6, J_{10}] \in \mathbb{WP}_{(2,4,6,10)}^3(\Ok)\setminus \{J_{10} = 0 \}$ be a given weighted moduli  point in the weighted projective space. We define the tuple $\p$ to be  a \textbf{normalized weighted moduli point} if the weighted greatest common divisor of its coordinates is one, where the weighted greatest common divisor  is defined as follows
\[  \wgcd (J_{2i} )  = \prod_{ v \in M_K}  \prod_{\substack{d^{q_i} | J_{i} \\ d \in \Ok}}  |d|_v \]
for all valuation $v$ in the set of valuations $M_K$ and divisors $d \in \Ok$ such that for all $i=2,4, 6, 10$ we have $d^i | J_{i}$.  Note that if $ K = \Q$ then,
\[  \wgcd (J_{2i} )  =    \prod_{\substack{p^{q_i} | J_{i} \\ p \in \Z}}  |p| \]
for all primes $p \in \Z$. The normalization of the point $\p \in \mathbb{WP}^3_\w(K)$ is
\[ \mathfrak q = \frac 1 {\wgcd (\p)} \star \p. \]
and from   \cite[Lemma.~1]{sh-h} we have that since $w= (2, 4, 6, 10)$   and $d = \gcd( 2, 4, 6, 10)=2$, then   this normalization is unique up to a multiplication by a $d$-th root of unity.  In other words, for every  normalized point  $\p= \left[J_2, J_4, J_6, J_{10}\right] \in \mathbb{WP}^3_\w(\Q)$, there is another normalized point  $\p^\prime = \left[-J_2, J_4, -J_6, -J_{10}\right]$ equivalent to $\p$. Moreover, $\p$ and $\p^\prime$ are isomorphic over the Gaussian integers, see \cite[Cor. 2]{sh-h} for more details.




\begin{rem}
Notice that if the first coordinate $J_2=0$ then the normalization  in the space $\mathbb{WP}_{(2,4,6,10)}^3 (\Q)$ is unique up to a multiplication by a $d=\gcd (4, 6, 10)$-th root of unity.  For example, the pair of points $\p$ and $\p^\prime$ in Table~\ref{tab-1} are twists of each other, where $\p^\prime = \sqrt{-1} \cdot \p$.   In the coming section we will see why these are the only tuples of height $h=1$ for which this is true.
\end{rem}

\begin{table}[htp]
\caption{Twists of height $\h=1$.}
\begin{center}
\begin{tabular}{|c|c|c|}
\toprule
\# & $\p$ & $\p^\prime =  \sqrt{-1} \cdot \p $ \\
\midrule
1 & [0, -1, -1, -1] & [0, -1, 1, 1] \\
2 & [0, -1, -1, 1] & [0, -1, 1, -1] \\
3 & [0, -1, 0, -1] & [0, -1, 0, 1] \\
4 & [0, 0, -1, -1] & [0, 0, 1, 1] \\
5 & [0, 0, -1, 1] & [0, 0, 1, -1] \\
6 & [0, 0, 0, -1] & [0, 0, 0, 1] \\
7 & [0, 1, -1, -1] & [0, 1, 1, 1] \\
8 & [0, 1, -1, 1] & [0, 1, 1, -1] \\
9 & [0, 1, 0, -1] & [0, 1, 0, 1] \\
\bottomrule
\end{tabular}
\end{center}
\label{tab-1}
\end{table}%

Next, we will see  when a point $ \p \in  \mathbb{WP}^3_w(\Q)$ is normalized over $\overline \Q$.  The  point $ \p$ is called \textbf{absolutely normalized}, i.e.  normalized over $\overline \Q$ if $\overline{\wgcd}(\p) = 1$ where
\[ \overline{\wgcd}(\p) =  \prod_{\substack{\l^{q_i} | J_{i} \\ \l \in \overline \Q}}  |\l|  \]
such that for all $i= 2, 4, 6, 10$, $\l^i \in \Z$. We can compute the absolutely normalization of a point in the weighted projective space the same way as we computed the normalization of that point.

Two absolutely normalized tuples are twists of each other if and only if they are isomorphic over the Gaussian integers.

The  \textbf{weighted multiplicative height} of $P=\left( x_0, \dots , x_n \right)  \in \P^n_w (K)$ over $K$ is defined as
\begin{equation}\label{height}
\wh_K (P)= \frac 1 { \wgcd(\p)} \,  \prod_{v \in M_K} \max  \left\{ |x_0|_v^{1/q_0},  \dots , |x_n|_v^{1/q_n}  \right\}
\end{equation}
for all valuations $v$ is the set of valuations   $M_K$ and the \textbf{absolute weighted multiplicative height}  over $\overline \Q$
\begin{equation}\label{abs-height}
\tilde \wh (P)= \frac 1 { \overline{\wgcd}(\p)} \,\prod_{v \in M_K}   \max  \left\{ |x_0|_v^{1/q_0},  \dots , |x_n|_v^{1/q_n}  \right\}.
\end{equation}

Note that if we consider points $\p \in \mathbb{WP}^3_w(\Q)$ then it is clear that $\p$  will have a representative over the integers.  With such a
representative for the coordinates of $\p$, the non-Archimedean absolute values give no contribution to the height and the only one that matters is the Archimedean absolute value, i.e. just the maximum of the coordinates of the given point.

Next, we show how we can calculate the absolute minimal tuple and absolute minimum height of all  zero-dimensional loci in $\M_2$. All such curves are defined over $\Q$, so we consider their respective weighted moduli points as existing in $\mathbb{WP}^3_\w(\Z)$. 
\begin{exa}
Let $\X$ be the genus two curve with equation $y^2 = x^6 - 1$. Its weighted moduli point is
\[
\begin{split}
\p_1 = \left[ 240, 1620, 119880, 46656 \right] = \left[ 2^4 \cdot 3 \cdot 5, 2^2 \cdot 3^4 \cdot 5, 2^3 \cdot 3^4 \cdot 5 \cdot 37, 2^6 \cdot 3^6 \right].
\end{split}
\]
which is a minimal tuple since there is no $p \in \Z$ such that $p^{2i} | J_{2i}$, for $i = 1, 2, 3, 5$. This tuple has three twists for values of $\l = 2, 3, 6$. The corresponding twists are
		\begin{alignat*}{2}
			\p_2 &=\frac 1 2 \star \p =  \left[ 120, 405, 14985, 1458 \right] &&= \left[ 2^3 \cdot 3 \cdot 5, 3^4 \cdot 5, 3^4 \cdot 5 \cdot 37, 2 \cdot 3^6 \right]
			\\
			\p_3 &= \frac 1 3 \star \p  = \left[ 80, 180, 4440, 192 \right] &&= \left[ 2^4 \cdot 5, 2^2 \cdot 3^2 \cdot 5, 2^3 \cdot 3 \cdot 5 \cdot 37, 2^6 \cdot 3 \right]
			\\
			\p_4 &= \frac 1 6 \star \p = \left[ 40, 45, 555, 6 \right] &&= \left[ 2^3 \cdot 5, 3^2 \cdot 5, 3 \cdot 5 \cdot 37, 2 \cdot 3 \right]
		\end{alignat*}
		\textit{The heights are}
		\[
			\h(\p_1) = 4\sqrt{15},\; \h(\p_2) = 2\sqrt{30},\; \h(\p_3) = 4\sqrt{5},\; \h(\p_4) = 2\sqrt{10}.
		\]
		The absolute minimal tuple is
		\[
			\begin{split}
				\p = \left[ 40, 45, 555, 6 \right] = \left[ 2^3 \cdot 5, 3^2 \cdot 5, 3 \cdot 5 \cdot 37, 2 \cdot 3 \right].
			\end{split}
		\]
		Hence, the absolute minimum height is
		\[
			\begin{split}
				\overline{\h}(\p) = 2\sqrt{10}.
			\end{split}
		\]
	\end{exa}
	
	\begin{exa}
		Let $\X$ be the genus two curve with equation $y^2 = x^5 - x$. Its weighted moduli point is
		\[
			\begin{split}
				\p_1 = \left[ -40, -80, 320, -256 \right] = \left[-1 \cdot 2^3 \cdot 5, -1 \cdot 2^4 \cdot 5, 2^6 \cdot 5, -1 \cdot 2^8 \right],
			\end{split}
		\]
		which is a minimal tuple since there is no $p \in \Z$ such that $p^{2i} | J_{2i}$, for $i = 1, 2, 3, 5$. This tuple has a sole twist for the value $\l = 2$. This twist is
		\[
			\begin{split}
				\p_2 = \left[ -20, -20, 40, -8 \right] = \left[-1 \cdot 2^2 \cdot 5, -1 \cdot 2^2 \cdot 5, 2^3 \cdot 5, -1 \cdot 2^3 \right].
			\end{split}
		\]
		The heights are
		\[
			\begin{split}
				\h(\p_1) = 2\sqrt{10},\; \h(\p_2) = 2\sqrt{5}.
			\end{split}
		\]
		The absolute minimal tuple is
		\[
			\begin{split}
				\p = \left[ -20, -20, 40, -8 \right] = \left[-1 \cdot 2^2 \cdot 5, -1 \cdot 2^2 \cdot 5, 2^3 \cdot 5, -1 \cdot 2^3 \right].
			\end{split}
		\]
		
		Hence, the absolute minimum height is
		\[
			\begin{split}
				\overline{\h}(\p) = 2\sqrt{5}.
			\end{split}
		\]
		
	\end{exa}
	
	\begin{exa}
		Let $\X$ be the genus two curve with equation $y^2 = x^6- x$. Its weighted moduli point is
		\[
			\begin{split}
				\p_1 = \left[ 0, 0, 0, 3125 \right] = \left[ 0, 0 ,0 , 5^5 \right].
			\end{split}
		\]
		which is a minimal tuple since there is no $p \in \Z$ such that $p^{2i} | J_{2i}$, for $i = 1, 2, 3, 5$. This tuple has a sole twist for the value $\l = 5$. This twist is
		\[
			\begin{split}
				\p_2 = \left[ 0, 0, 0, 1\right]
			\end{split}
		\]
		The heights are
		\[
			\begin{split}
				\h(\p_1) = \sqrt[\leftroot{-2}\uproot{2}10]{3125},\; \h(\p_2) = 1
			\end{split}
		\]
		The absolute minimal tuple is
		\[
			\begin{split}
				\p = \left[ 0, 0, 0, 1\right].
			\end{split}
		\]
		Hence, the absolute minimum height is
		\[
			\begin{split}
				\overline{\h}(\p) = 1.
			\end{split}
		\]
	\end{exa}

In \cite{height} was defined the moduli  height  of genus $g$ algebraic curves.  Consider  $ \p \in \M_g$, where $\M_g$ is the moduli space of smooth  irreducible algebraic curves of genus $g \geq 2$.  It is known that $\M_g$ is a quasi projective variety of dimension $3g-3$ and is embedded in $\P^{3g-2}$. The moduli height is  the usual height in the projective space $\P^{3g-2}$, i.e. if $P \in \P^{3g-2} (\Q) $ then
\[ H(P)  = \max \left\{ \frac{}{}  |x_0|,  \dots , |x_n|  \right\}.  \]
Moreover, for  $g = 2$  where a moduli point is given by  $\p = [J_2, \dots, J_{10}] \in \M_2(\Q)$ then
\[ H(\p)  =  \max \left\{  \frac{}{} |J_2|,  \dots , |J_{10}|  \right\}. \]
If we consider $\p \in \mathbb{WP}_w^3 (\Q)$ then the weighted moduli height is
\[ \h (\p) = \max \left \{ |J_2|^{1/2}, \dots, |J_{10}|^{1/10}    \right  \}.\]

Hence, by definition for a given moduli point  the weighted moduli height is always smaller than the moduli height.  And this makes computing databases of genus two curves based on the weighted moduli height much more convenient.

\subsection{Naive height versus weighted moduli height}

Let $\X$ be a genus two curve with equation given as
\[ y^2 =  a_6 x^6 + a_5x^5z + \cdots + a_0 z^6.    \]
The naive height $H(\X)$ of a curve  (resp. polynomial) is the maximum of the absolute value of the  coefficients of the curve (resp. polynomial);  see \cite{height} for details.  In \cite{height} and \cite{data} were used both the naive height and the moduli height to analyze points in $\M_2 (\Q)$.  Even though   using the weighted moduli height makes more sense in analyzing points in $\M_2 (\Q)$,  since a weighted moduli point determines uniquely an isomorphism class of the curves, it would be interesting to see how the weighted moduli height relates to the naive height.    In \cite{height} it  is given an inequality how the naive height is compared to the  moduli height. In this subsection we capture previous results and also compare the naive height to the weighted moduli height.

 In \cite{height}  the naive height of a curve is compared with the moduli height and  showed that  for a genus 2 curve with equation $y^2=f(x)$ the moduli height is bounded as follows
\[ \tilde \h (f) \leq 2^{28} \cdot 3^9 \cdot 5^5 \cdot 7 \cdot 11 \cdot 13 \cdot  17 \cdot 43  \cdot H(f)^{10}. \]
Note that there is a typo in the computations in \cite{height} and the above coefficient should be $   2^{24}\cdot 3^7 \cdot 5^4 \cdot 7 \cdot 11 \cdot 13 \cdot 43$.  Next,  we consider how the weighted moduli height is related to the naive height of the curve. 

\begin{lem}\label{bound-gen-2}
Let $\X$ be a genus 2 curve and $\p$ its corresponding point in $ \mathbb{WP}_w^3 (\Q)$.  Then   the weighted moduli height is bounded as follows
\[  \h (\p) \leq 2^3  \sqrt {3 \cdot 5 \cdot 7}  \cdot H (\X)\]
\end{lem}
\proof We assume that the equation of the curve is given as $y^2 = f(x)$. Recall  that each $J_{2i}$ is a degree $2i$ polynomial evaluated at $f$, i.e degree $2i$ polynomial  given in   $k[a_0, \cdots, a_6]$.  Then, from  Lemma 15 in  \cite{height} we have
\[H    (J_{2i}      (f) )\leq c_0 \cdot H    (J_{2i})  \cdot H(f)^{2i}  \]
where $c_0 =  {d+n \choose n } $ is the  number of monomials of a degree $d$  homogenous polynomial in $n+1$ variables.  Computations of $H(J_{2i})$ is done in Maple and we get
\[
\begin{split}
H(J_2) & \leq  2^{6} \cdot 3  \cdot 5 \cdot 7 \cdot H(f)^2 \\
H(J_4) & \leq 2^{3} \cdot 3^{5} \cdot 5^2 \cdot 7 \cdot H(f)^4  \\
H(J_6) & \leq  2^5 \cdot 3^5 \cdot 5  \cdot 7 \cdot 11 \cdot 37 \cdot H(f)^6   \\
H(J_{10}) & \leq 2^{9} \cdot 3^5 \cdot 5  \cdot 7 \cdot 11 \cdot 13 \cdot H(f)^{10}  \\
\end{split}
\]
The weighted moduli height of $\p$ is computed as follows
\[ \h (\p) = \max \{ H (J_2(f))^{1/2}, \dots , H( J_{10} (f))^{1/10}\} \leq  2^3 \sqrt {3 \cdot 5 \cdot 7} \cdot H(f).  \]
This concludes the proof.
\qed

Hence,  we get that the height of the given curve  with corresponding moduli point $\p =[J_2: J_4: J_6: J_{10}]$ is as follows
\[ H (\X)  \geq \frac{1}{2^3  \sqrt {3 \cdot 5 \cdot 7}}   \cdot \h (\p)  \]
 but  we don't know if this lower bound is strict and if it can ever be achieved.   In \cite{height} are displayed all genus 2 curves (up to isomorphism) with naive height $H =1$. In \cite[Tables 1-4]{height} there are displayed 230 isomorphism classes of genus two curves with naive height one.  There are 184 genus two curves with naive height one and automorphism group $\Z_2$, 30 with automorphism group isomorphic to $V_4$, 11 with automorphism group $D_4$, two curves with group $D_6$, and three with automorphism group respectively 10, 24, and 48.

 In an upcoming paper we will  check how the weighted moduli height of those curves in \cite[Tables 1-4]{height}  is compared to their naive height and see if this lower bound is sharp for any of those curves.  

\section{Rational points in the weighted moduli space}\label{sec:5}

Next we will build a database of all isomorphism classes of genus 2 curves with bounded weighted moduli height.  We  start with creating a database of normalized tuples  in $\mathbb{WP}_w^3 (\Q)\setminus \{J_{10} =0\}   $ of height $\leq h$, for a given $h$.    Then,  we can compute the equation of the curve for a given $\p \in \mathbb{WP}_w^3 (\Q)\setminus \{J_{10} =0\} $ using methods explained in the upcoming  subsection~\ref{sec:4}.



Let $h$ be a positive integer.  The number of  points $\p \in \WP (\Q) \setminus \{J_{10} = 0\}$ of  weighted moduli height $\leq h$ is given by the following lemma.
	
\begin{lem} For any given positive integer $h$ there are  at most
\[  2 h^{10} (h^2 + 1) (2h^4 +1) (2h^6 +1) \]	
normalized points in  $  \WP (\Q) \setminus \{J_{10} = 0\}$   of  weighted moduli height $\leq h$.
\end{lem}	

\proof   Let $\X$ be a genus two curve and   $\p = [J_2, J_4, J_6, J_{10}] \in \WP (\Q) \setminus \{J_{10} = 0\}$ the moduli point representing $\X$ of height $\h(\p) = h$.   From the definition of weighted moduli height  for every $i=1, 2, 3, 5$ we have that
\[ |J_{2i} | ^{1/2i} \leq h\]
Hence, for a given height $h$, every element $J_{2i}$ of $\p$ satisfies
 \[ -h^{2i} \leq J_{2i} \leq h^{2i}  \text{ for  } i = 1, 2, 3, 5.\]
Since,   $J_{10}\neq 0$ then there are $2h^{10}$ possible integer values for $J_{10}$.  In \cite[Corr.2]{sh-h} it is proved that for every  normalized point  $\p= \left[J_2, J_4, J_6, J_{10}\right] \in \mathbb{WP}_{(2, 4, 6, 10)}^3 (\Q)$, there is another normalized point  $\p^\prime = \left[-J_2, J_4, -J_6, -J_{10}\right]$ equivalent to $\p$.  Since, we only want to count equivalence classes then we let $J_2 \geq 0$ and therefore there are $h^2+1$ possible integer values for $J_2$.  Lastly, for  $i = 2, 3$ there are $(2 h^{2i} +1)$-possible integer values for $J_{2i}$.

But, up to this point we didn't count the points $\p \in \mathbb{WP}_{(2, 4, 6, 10)}^3 (\Q)$ that satisfy the property that there exists a $\lambda \in \Q^\star$ such that
\begin{equation}\label{lambda-star}
\lambda \star (J_2,  \dots, J_{10})  =  (\lambda^2 J_2, \dots, \lambda^{10} J_{10}).
\end{equation}

But there is not a rigorous mathematical way of counting this points. Hence, the number of normalized points $\p \in \mathbb{WP}_{(2, 4, 6, 10)}^3 (\Q)$ will be at most
\[  2 h^{10} (h^2 + 1) (2h^4 +1) (2h^6 +1). \]	
 This concludes the proof.
\qed

 However, for computational reasons it is more convenient to consider the space $\mathbb{WP}_{(1,2,3,5)}^3 (\Q)$ instead of $\WP (\Q)$. In this case the above lemma becomes as follows.

\begin{lem}\label{lemma-1235} For any given positive integer $h$ there are  at most
\[  2 h^{5} (h + 1) (2h^2 +1) (2h^3 +1) \]	
normalized points in  $\mathbb{WP}_{(1,2,3,5)}^3 (\Q)$   of  weighted moduli height $\leq h$.
\end{lem}


\begin{table}[hbp]
\centering
\caption{Points of bounded height}
\label{table-3}
\begin{tabular}{|c|c|c|}
\toprule
height & \# of tuples & \# of points in $\mathbb{WP}_{(1, 2, 3, 5)}^3 (\Q)$ \\
\midrule
1 & 36 & 27\\
2 & 29 376 & 24 423\\
3 & 2 031 480 & 1 776 549\\
4 & 43 591 680 & 39 206 865\\
\bottomrule
\end{tabular}
\end{table}

The proof is similar to the previous lemma.  In Table~\ref{table-3}  for  weighted moduli height $\h \leq 4$ we display how many possible tuples are there in $\mathbb{WP}_{(1, 2, 3, 5)}^3 (\Q)$  and then  in the next column we display how many of this tuples give us isomorphism classes of genus two curves. And this results agree with Lemma~\ref{lemma-1235}.

Then, in the upcoming  Table~\ref{table-0}  we display  all points $\p \in  \mathbb{WP}_{(1,2,3,5)}^3 (\Q) \setminus  \{J_{10} = 0\} $ with weighted moduli heights $\h = 1$.

\begin{small}
\begin{table}[h!]
\centering
\caption{Tuples with weighted moduli height $\h=1$}
\begin{tabular}{|c|c|c|c|c|c|}
	\toprule
	\# & $\p=[J_2: J_4: J_6: J_{10}]$    & \#   &    $\p=[J_2: J_4: J_6: J_{10}]$  & \#   &    $\p=[J_2: J_4: J_6: J_{10}]$  \\
	\midrule
	1 & [0, -1, 0, 1]  & 	10 &  [1, 0, 1, 1]  &	19 & [0, 1, 1, 1] \\
	2 & [0, 1, 0, 1]    & 	11 & [1, -1, -1, 1] & 	20 & [1, 0, 1, -1]\\
	3 & [0, -1, 1, 1] & 12 &  [1, 1, -1, 1] & 21 & [1, -1, -1, -1]  \\
	4 & [0, 0, 0, 1] & 13 &  [1, 1, 1, -1] & 22 & [1, 1, -1, -1]\\
	5 & [0, 0, 1, -1] & 14 & [1, -1, 1, -1] & 23 & [1, -1, 0, -1]\\
	6 & [0, 0, 1, 1] & 15 & [1, 1, 1, 1] & 24 & [1, 1, 0, -1]\\
	7 & [1, 0, -1, 1]   & 16 & [1, 0, -1, -1] & 25 & [1, 1, 0, 1]\\
	8 & [1, 0, 0, -1]  & 17 & [0, -1, 1, -1] & 	26 & [1, -1, 0, 1] \\
	9 & [1, 0, 0, 1] & 18 & [0, 1, 1, -1] & 	27 & [1, -1, 1, 1] \\
%
	\bottomrule
\end{tabular}
\label{table-0}
\end{table}
\end{small}

Note that  in Table~\ref{table-0}  the points $1$ to $15$ are fine moduli points, i.e. with methods that we will explain in the upcoming subsection we can compute their equations defined over $\Q$.    We have already seen points $1$ to $6$ and $17$ to $19$ and their twists in the previous section in Table~\ref{tab-1}.
	
\subsection{Equations of curves from points in $\WP (\Q)$ }\label{sec:4}

In \cite{m-sh} for every point $\p =[J_2 : J_4 : J_6 :J_{10}] \in 	\WP \setminus \{J_{10}=0\}$ was constructed a genus two curve defined over a field of definition. We recapture that result below:

		Let $\p = [J_2, J_4, J_6, J_{10}] \in \mathbb{WP}^3_\w(\Ok)$ for some number field $K$ and $\Ok$ its ring of integers. There is a genus-two curve corresponding to $\p$ as follows:
		\begin{enumerate}[i)]
			\item If $J_2 \cdot J_{10} \neq 0$, then there is a genus 2 curve $C_{(\alpha, \beta)}$ given by
			\[ C_{(\alpha, \beta)} : y^2 = \sum_{i=0}^6 a_i (\alpha, \beta) x^i\]
			
			with coefficients given in  \cite[Eq. 3.9 and Eq. 3.10]{m-sh} and a pair $(\alpha, \beta)$ satisfying
			\[
				\alpha^2 + \l_6 \beta^2 \sigma = \gamma,
			\]
			where $\l_6$, $\sigma$, and $\gamma$ are given in terms of some $\rho$ and $\kappa$, which are themselves determined by $\p$. Their equations are given in \cite{m-sh}.  Moreover, $C_{(\alpha, \beta)}$ is defined over its field of moduli $K$, i.e. $a_i(\alpha, \beta) \in K$, $i = 0, \dots, 6$, if and only if $K$-rational $\alpha$ and $\beta$ exist. For the explicit equation of $C_{(\alpha, \beta)}$ we refer the reader directly to \cite[Theorem~1]{m-sh}.
			
			\item If $J_2 = 0$ and $J_4 \cdot J_6 \cdot J_{10} \neq 0$, then there is a genus 2 curve given by setting $\rho = \kappa \neq 0$ in the previous equation.
			
			\item If $J_2 = J_6 = 0$ and $J_4 \cdot J_{10} \neq 0$, there is only one genus-two curve given by
			\[
				\begin{split}
					y^2 =& (4\nu + 1)(2\nu - 1)x^6 + 2(1 - \nu)(4\nu + 3)x^5 - 15(1 - \nu)x^4 \\
					&+ 20(1 - \nu)^2x^3 + 5(2\nu - 3)(1- \nu)^2x^2 + 6(1 - \nu)^3x - (1 - \nu)^3,
				\end{split}
			\]
			where
			\[
				\nu = \frac{J_4^5}{2^23^55^5J_{10}^2}.
			\]
			
			\item  If $J_2 = J_4 = 0$ and $J_6 \cdot J_{10} \neq 0$, there is only one genus-two curve given by
			\[
				\begin{split}
					y^2 =& 5x^6 + 12(1 - \mu)x^5 - 15(1 - \mu)x^4 - 80(1 - \mu)^2x^3 \\
					&+ 15(4\mu - 7)(1 - \mu)^2x^2 - 60(1 - \mu)^3x + (4\mu - 13)(1 - \mu)^3,
				\end{split}
			\]
			where
			\[
				\mu = \frac{J_6^5}{2^43^45^5J_{10}^3} .
			\]
			
			\item If $J_2 = J_4 = J_6 = 0$ and $J_{10} \neq 0$, there is only one genus-two curve given by
			\[
				y^2 = x^6 - x.
			\]
		\end{enumerate}

Now that we have seen how to compute equations of curves  we turn back to one of our main goals for this paper: Investigate what proportion of points in $\WP (\Q)$ are defined over $\Q$.  The following remark is believed to be true.

\begin{rem}\label{thm-16}
For almost all  points  $\p \in \M_2 (\Q)$,  the field of moduli is not a field of definition. In other words, the probability that a point $\p \in \M_2 (\Q)$ has a representative curve defined over $\Q$ is zero.
\end{rem}

 However, this was contradicted by previous results in \cite{data} where absolute invariants $i_1, i_2, i_3$ were used to identify the isomorphism classes of curves. At first we believed, this  was because such invariants are not everywhere defined and to build the database two different set of invariants were used, see \cite{data} for more details.    But,  as we will see below  we get the same contradicting results using points in $\mathbb{WP}_{(1, 2, 3, 5)}^3 (\Q)$ with minimal weighted moduli height.
The following result  is well known by work of Cardona/Quer \cite{cardona} or  Shaska; see \cite{ants-5}.

\begin{lem}
Every genus 2 curve with automorphism group of order $> 2$ is defined over its field of moduli. Moreover, its equation over its field of moduli
it's  given in \cite[Lemma~6]{g-sh}.
\end{lem}

\begin{table}[hbp]
\caption{The number of points in $\mathbb{WP}_{(1, 2, 3, 5)}^3 (\Q)$ for each automorphism group. }
\begin{center}
\begin{tabular}{|c|c|c|c|c|c|c|c|c|}
\toprule
$\awh$  &  \# pts & $D_4$ & $D_6$ &  $V_4$ & $C_2$ & $C_{10}$ & \#fine pts.  &  $\frac{\#\text{fine pts.}}{\#\text{pts.} }$\\
\midrule
$0 < \h \leq 1$  & 27 & 0 & 0 & 0 & 26 & 1 & 15 & 0.555\\
$1 < \h \leq 2$  & 24 396 & 0 & 0 & 0 & 24 396 & 0& 9 423 & 0.386 \\
$2 < \h \leq 3$ & 1 752 126 & 0 & 0 & 2 & 1 752 124 & 0  & 596 818  &  0.340 \\
$3 < \h \leq 4$  & 37 430 316 & 1 & 0 & 14 & 37 430 301 & 0 & 11 365 256 & 0.303  \\
\bottomrule
\end{tabular}
\end{center}
\label{default}
\end{table}%

Hence, we are left to check for curves with automorphism group of order two.
In  Table~\ref{default} given below we  display the number of points with height $\awh \leq 5$ for each possible group of automorphisms and  the number of points that  are defined over the rationals.  We stop at $\h =5$ since then the computations get to big.  As we can see most of the curves are curves with automorphism group $\Z_2$, as expected.   From the data we notice that the  percentage of the curves defined over the rationals  gets smaller as the weighted moduli height gets higher  and as we can see  approximately $30\%$ of the points are defined over the rationals.  This results agree with previous results from \cite{data} and we expect the ratio to converges to zero as the height gets bigger. It would be interesting  to see   for what moduli height $\h$ this ratio is zero.



\bigskip

\noindent \textbf{Acknowledgments.}  We want to thank Prof. Tony Shaska for sharing his code for the universal equation of genus two curves from \cite{m-sh}  and his helpful discussions and suggestions.


\bibliographystyle{amsplain}
	
\bibliography{ref}{}

\begin{bibdiv}
\begin{biblist}

\bib{data}{article}{
      author={Beshaj, L.},
      author={Hidalgo, R.},
      author={Malmendier, A.},
      author={Kruk, S.},
      author={Quispe, S.},
      author={Shaska, T.},
       title={Rational points on the moduli space of genus two},
        date={2018},
     journal={Contemporary Math.},
      volume={703},
       pages={87\ndash 120},
}

\bib{cardona}{incollection}{
      author={Cardona, Gabriel},
      author={Quer, Jordi},
       title={Field of moduli and field of definition for curves of genus 2},
        date={2005},
   booktitle={Computational aspects of algebraic curves},
      series={Lecture Notes Ser. Comput.},
      volume={13},
   publisher={World Sci. Publ., Hackensack, NJ},
       pages={71\ndash 83},
         url={https://doi.org/10.1142/9789812701640_0006},
      review={\MR{2181874}},
}

\bib{g-sh}{article}{
      author={Gutierrez, J.},
      author={Shaska, T.},
       title={Hyperelliptic curves with extra involutions},
        date={2005},
        ISSN={1461-1570},
     journal={LMS J. Comput. Math.},
      volume={8},
       pages={102\ndash 115},
         url={https://doi.org/10.1112/S1461157000000917},
      review={\MR{2135032}},
}

\bib{Ig-60}{article}{
      author={Igusa, Jun-Ichi},
       title={Arithmetic variety of moduli for genus two},
        date={1960},
        ISSN={0003-486X},
     journal={Ann. of Math. (2)},
      volume={72},
       pages={612\ndash 649},
      review={\MR{0114819}},
}

\bib{m-sh-2}{article}{
      author={Malmendier, A.},
      author={Shaska, T.},
       title={The {S}atake sextic in {F}-theory},
        date={2017},
        ISSN={0393-0440},
     journal={J. Geom. Phys.},
      volume={120},
       pages={290\ndash 305},
         url={https://doi.org/10.1016/j.geomphys.2017.06.010},
      review={\MR{3712162}},
}

\bib{m-sh}{article}{
      author={Malmendier, Andreas},
      author={Shaska, Tony},
       title={A universal genus-two curve from {S}iegel modular forms},
        date={2017},
        ISSN={1815-0659},
     journal={SIGMA Symmetry Integrability Geom. Methods Appl.},
      volume={13},
       pages={089, 17 pages},
         url={https://doi.org/10.3842/SIGMA.2017.089},
      review={\MR{3731039}},
}

\bib{Me}{incollection}{
      author={Mestre, Jean-Fran{\c{c}}ois},
       title={Construction de courbes de genre {$2$} {\`a} partir de leurs
  modules},
        date={1991},
   booktitle={Effective methods in algebraic geometry ({C}astiglioncello,
  1990)},
      series={Progr. Math.},
      volume={94},
   publisher={Birkh{\"a}user Boston, Boston, MA},
       pages={313\ndash 334},
      review={\MR{1106431}},
}

\bib{sh-h}{article}{
      author={Shaska, T.},
       title={Heights on weighted projective spaces},
        date={201801},
     journal={submitted},
       pages={1\ndash 10},
      eprint={1801.06250},
         url={https://arxiv.org/abs/1801.06250},
}

\bib{height}{incollection}{
      author={Shaska, T.},
      author={Beshaj, L.},
       title={Heights on algebraic curves},
        date={2015},
   booktitle={Advances on superelliptic curves and their applications},
      series={NATO Sci. Peace Secur. Ser. D Inf. Commun. Secur.},
      volume={41},
   publisher={IOS, Amsterdam},
       pages={137\ndash 175},
      review={\MR{3525576}},
}

\bib{deg2}{article}{
      author={Shaska, T.},
      author={V\"olklein, H.},
       title={Elliptic subfields and automorphisms of genus 2 function fields},
        date={2000},
     journal={Algebra, arithmetic and geometry with applications},
       pages={703\ndash 723},
}

\bib{ants-5}{incollection}{
      author={Shaska, Tony},
       title={Genus 2 curves with {$(3,3)$}-split {J}acobian and large
  automorphism group},
        date={2002},
   booktitle={Algorithmic number theory ({S}ydney, 2002)},
      series={Lecture Notes in Comput. Sci.},
      volume={2369},
   publisher={Springer, Berlin},
       pages={205\ndash 218},
         url={https://doi.org/10.1007/3-540-45455-1_17},
      review={\MR{2041085}},
}

\end{biblist}
\end{bibdiv}
	
\end{document}